\documentclass[aap]{imsart}

\RequirePackage{amsthm,amsmath,amsfonts,amssymb}
\RequirePackage[numbers]{natbib}
\RequirePackage[colorlinks,citecolor=blue,urlcolor=blue]{hyperref}
\RequirePackage{graphicx}

\usepackage{stmaryrd}
\usepackage{lmodern}

\startlocaldefs

\newtheorem{theorem}{Theorem}[section]
\newtheorem{lemma}[theorem]{Lemma}
\newtheorem{proposition}[theorem]{Proposition}
\newtheorem{corollary}[theorem]{Corollary}
\theoremstyle{remark}
\newtheorem{definition}[theorem]{Definition}
\newtheorem{remark}{Remark}

\newcommand{\E}        {{ {\rm I \hskip -2pt E}}}
\newcommand{\N}        {{ {\rm I \hskip -2pt N}}}

\newcommand{\X}{\mathcal{X}}
\newcommand{\Z}        {{\mathbb Z}}
\renewcommand{\P}      {{ {\rm I \hskip -2pt P}}}
\newcommand{\R}        {{\rm I \hskip -2pt R}}
\newcommand{\pib}   {\overline{\pi}}
\newcommand{\pibp}  {\overline{\pi}^+}
\newcommand{\M}  {{\mathcal M}}
\newcommand{\MM} {\overline{\mathcal M}}

\newcommand{\indic}[1]{1_{\left\lbrace #1 \right\rbrace}}
\endlocaldefs

\begin{document}

\begin{frontmatter}
\title{Household epidemic models and McKean-Vlasov Poisson driven stochastic differential equations}
\runtitle{Household epidemic models}

\begin{aug}
\author[A]{\fnms{Rapha\"el} \snm{Forien}\ead[label=e1]{raphael.forien@inrae.fr}}
\and
\author[B]{\fnms{\'Etienne} \snm{Pardoux}\ead[label=e2]{etienne.pardoux@univ-amu.fr}}
\address[A]{INRAE, Centre INRAE PACA, Domaine St-Paul - Site Agroparc 84914 Avignon Cedex France,
\printead{e1}}

\address[B]{Aix–Marseille Université, CNRS, I2M, UMR 7373 13453 Marseille, France,
\printead{e2}}
\end{aug}

\begin{abstract}
This paper presents a new view of household epidemic models, where we exploit the fact that the interaction between the households is of mean field type. We prove the convergence, as the number of households tends to infinity, of the number of infectious individuals in a uniformly chosen household to a nonlinear Markov process solving a McKean--Vlasov Poisson driven stochastic differential equation, as well as a propagation of chaos result. We also define a basic reproduction number $R_\ast$ and show that if $R_\ast>1$, then the nonlinear Markov process has a unique non trivial ergodic invariant probability measure, whereas if $R_\ast\le1$, it converges to $0$ as $t$ tends to infinity. 
\end{abstract}



\end{frontmatter}

\section{Introduction} \label{sec:intro}
In this paper, we present a new view of household epidemic models. 
Motivated by its simplicity, we present it in the particular case of the SIS model, but the same approach can be developed for other types of epidemic models, like the SIR, SIRS, SIR model with demography, and others. 
We recall that S stands for susceptible, I for infectious and R for removed. 
The SIS model describes an epidemic where individuals who become infected (i.e. move from the compartment of susceptibles to the compartment of infected/infectious individuals) become susceptible again when they recover. 
The SIR describes an epidemic where an infectious individual who recovers becomes removed, which means that it is immune and cannot be infected ever again, while in the SIRS model, removed individuals lose their immunity after some time, and eventually become susceptible again. 
In all those models, the total population size remains constant, while in the ``SIR model with demography'', the total population size fluctuates, due to births and deaths. 
Note that in the SIS, SIRS models and in the SIR model with demography, there is constantly a flux of new susceptibles, which allows the disease to become endemic, unlike in the SIR case, where the total number of individuals infected during an epidemic cannot exceed the number of individuals who are susceptible at the start of the epidemic.

A very natural step in changing homogeneous epidemic models into more realistic models is to include households, which are small groups of individuals who interact more frequently within their group than with other individuals in the population.
Household models are thus a key example of two--level mixing models. 
This describes both the situation of human populations, but also of many domestic animal populations, where cages/sheds in poultry farms or pens/fields in sheep/cattle farms play the role of households.

Household models can be roughly described as follows. The total population is the (disjoint) union of households of relatively small (and varying) size. 
Each infectious individual infects any other individual in the same household at a ``local rate'' $\lambda_L$, and any other individual in the total population at a ``global rate'' $\lambda_G$. 
In the last sentence, ``any other'' means ``chosen uniformly at random''. 
The infectious periods are i.i.d., in our case exponential with a given parameter $\gamma$ (so that the model is Markov). 
Note that both local and global infections are also often considered in the case of epidemics on large graphs, see e.g. Neal \cite{Neal2006} who studies an SIS model on a large circle with local and global infectious contacts and establishes limit theorems for the total proportion of infected individuals as the size of the circle tends to infinity.

The first papers on epidemic models with two levels of mixing go back to the 1950's, with Rushton and Mautner \cite{RM}
who study deterministic models, Bartlett \cite{Bar} and Daley \cite{Da} who study  stochastic models. We refer to 
Ball, Mollison and Scalia--Tomba \cite{BMS}, Ball \cite{Ball}, Neal \cite{Neal2006}, Ball, Britton and Neal \cite{BBN} who give a deep study of stochastic SIR and SIS epidemic models with two levels of mixing, as well as to Ball and Sirl \cite{BS} for an up--to--date presentation of stochastic SIR epidemics in structured populations, and for more references.  

Our viewpoint in this paper is to study asymptotic results as the number of households (and hence also the total population size) tends to infinity, while the household sizes remain unchanged. 
It is easy to see that the interactions between the various households is of {\it mean--field} type. 
This is reminiscent of the situation of particle systems which were studied by Sznitman \cite{ASS}. 
We establish a result of {\it propagation of chaos}, and prove that in the limit of an infinite number of households, the typical epidemic in a household follows a so--called nonlinear Markov process, whose transitions depend not only upon the situation of the epidemic in the household, but also upon its probability law through its mean, which is the limiting effect of the infections coming from the other households. 
Similar non--linear Markov processes have a long history, with in particular the work of McKean \cite{McK}. 
The SDEs of those nonlinear Markov processes are called McKean--Vlasov SDEs. 

Most of the literature on that topic treats Brownian--driven SDEs. 
However, L\'eonard \cite{Le} considers an epidemic model where the infection is the effect of a mean field interaction, and he obtains a McKean--Vlasov type SDE of Poissonian type as a law of large numbers limit.
Non-linear Markov processes with jumps have also been studied, notably in \cite{Gra1992a,Gra1992b} under uniform Lipschitz conditions on the jump rates and in \cite{Mel1996} with a specific form of interaction, see also \cite{APF}.
In our system, however, the jump rates are unbounded and not globally Lipschitz because household sizes may be arbitrarily large (see equation~\eqref{Nsde} below).
More recently, Carmona and Wang \cite{CW1,CW2} have studied finite state mean--field games and applications to contract theory and epidemics. 
Note that mean--field games have already been applied to malware epidemics by Bensoussan and Kolokoltsov \cite{KB}, see also M.H.R. Khouzani, S. Sarkar and E. Altman \cite{KSA} on a similar topic.
Related models have also been studied in the field of neuroscience, as for example in \cite{de_masi_hydrodynamic_2015,fournier_toy_2016,chevallier_mean-field_2017}.

Should we assume that the household sizes were bounded, then the existence and uniqueness of the nonlinear Markov process would be very elementary. Indeed, the Fokker--Planck equation for the evolution of its law would be a finite dimensional system of ODEs with locally Lipschitz coefficients, whose solution cannot explode since it is a probability distribution. Once all the marginal laws of the process are specified, then the SDE for the nonlinear Markov process becomes a classical easy to solve Poisson driven SDE. However, we only assume that the household size is a square integrable random variable. We prove existence and uniqueness of a solution of the nonlinear Fokker--Planck equation by a fixed point argument and approximation by both a decreasing sequence of supersolutions, and an increasing sequence of subsolutions. 
Note that allowing rare arbitrarily large households  can account for situations where some households are much larger than typical households, as is the case for example when considering schools or workplaces.

We next study the large time behavior of our limiting SDE. We define the basic household reproduction number $R_\ast$, which is the mean number of households infected as a result of a local infection in a typical infected household, started with one infectious individual.
We give an explicit formula for $R_\ast$.
If $R_\ast\le1$, then the number of infectious individuals in a typical household tends to 0 as $t\to\infty$, whereas if $R_\ast>1$, the law of that number converges to an invariant measure which is not the Dirac measure at zero. 

Let us outline our strategy for proving this result. Our nonlinear Markov process $(X(t), t \geq 0)$ is the solution of an SDE whose coefficients depend on the law of $X(t)$ through its mean $\E[X(t)]$.
The long time behaviour of these non--linear Markov processes is much harder to study than for classical Markov processes, as the usual ergodic theory does not apply.
To circumvent this difficulty, we first consider a so--called ``forced process'' $(X_t(m), t \geq 0)$, where $\E[X(t)]$ is replaced by a fixed function $m(t)$, and show that the limit as $t\to\infty$ of the law of $X_t(m)$ depends only on $m(\infty):=\lim_{t\to\infty}m(t)$. 
We then prove that the limit as $t\to\infty$ of the law of $X(t)$ is the law of $X_\infty(m_\star)$, where $m_\star$ solves the fixed point problem $m_\star=\E[X_\infty(m_\star)]$ and that this fixed point must be $0$ if $R_\ast\le1$, and that there exists just one other fixed point $m_\star>0$ if $R_\ast>1$. 
We then show that if $\E[X(0)]>0$ and $R_\ast>1$, then the process does indeed converge in distribution towards $ X_\infty(m_\star) $. 
For that purpose, we use a minoration by a supercritical non--Markov branching process, to show that the proportion of infected households cannot approach zero.

Our results extend well--known classical results concerning the case where all households have size $1$ (the homogeneous model). Note that we shall study the fluctuations around the law of large numbers obtained here in another publication.

The paper is organized as follows. The model is defined precisely in Section~2. Section~3 states the three main results of the paper, namely Theorem \ref{thm:existence_uniqueness} which gives the existence and uniqueness of the nonlinear Markov process, Theorem~\ref{thm:chaos} which states the propagation of chaos result (which might be considered as a law of large numbers), and finally Theorem~\ref{thm:large_time} which gives the large time behavior of the nonlinear Markov process. Section 4 studies what we call the ``forced process'', which is our nonlinear Markov process, where we replace the unknown quantity $\E[X(t)]$ by a given function $m(t)$. In particular, we establish the monotonicity property of the forced process as a function of $m$. That property is exploited in an essential way in Section~5 for the proof of Theorem~\ref{thm:existence_uniqueness}. Section~6 is devoted to the proof of Theorem~\ref{thm:chaos} and finally Section~7 to the definition and computation of $R_\ast$, and the proof of Theorem~\ref{thm:large_time}.  In this last section, we use in particular  a comparison with a non--Markov continuous time branching process.

\section{Definition of the model} \label{sec:model}

We consider an SIS household epidemic model. 
In our model, the population consists of $N$ households, with sizes $\nu_1,\nu_2,\ldots,\nu_N$, where the $\nu_i$'s are i.i.d. $\N$--valued random variables.
Let $X_i^N(t)$ denote the number of infectious individuals in the $i$--th household at time $t$.

We suppose that each infectious individual ``tries'' to infect another individual within the same household at rate $ \lambda_L $, for some $ \lambda_L > 0 $ (the infected individual is chosen uniformly from those in the household). If this individual is susceptible, it becomes infected, and if it is already infectious, nothing happens. 
Moreover, each infectious individual ``tries'' to infect another individual chosen \textit{uniformly from the whole population} at rate $ \lambda_G $, for some $ \lambda_G $ (again, if it is already infectious nothing happens).
Finally, each infectious individual becomes susceptible at rate $ \gamma $, for some $ \gamma > 0 $.
The parameters $ \lambda_L $ and $ \lambda_G $ are the rates of local (respectively global) infections. 

Given the state of the process $ ((\nu_i, X^N_i(t)), 1 \leq i \leq N) $ at some time $ t \geq 0 $, a local infection takes place in the $ i $-th household at instantaneous rate
\begin{align*}
\frac{\nu_i-X^N_i(t^-)}{\nu_i-1} \lambda_L X^N_i(t^-),
\end{align*}
where the first factor corresponds to the probability that the chosen individual is susceptible (note that if $ \nu_i = 1 $, then $ (\nu_i-X^N_i(t))X^N_i(t) = 0 $, and we take the convention that the above rate is zero).
We then note that, for each global infection, choosing an individual uniformly from the population is equivalent to first choosing a household from the size-biased distribution and then choosing an individual uniformly in this household.
Hence the probability that a global infection targets an individual in the $ i $-th household is
\begin{align*}
\frac{\nu_i}{\sum_{j=1}^{N} \nu_j}.
\end{align*}
Global infections take place at rate $ \lambda_G \sum_{i=1}^{N} X^N_i(t^-) $, and, given that a global infection targets the $ i $-th household, it results in an actual infection with probability $ \frac{\nu_i-X^N_i(t^-)}{\nu_i} $ (for simplicity, we assume that an individual may try to infect itself during a global infection, contrary to the way we model local infections, noting that this happens with probability $ 1/N $).

Both for local and global infections, our model assumes that each infectious individual meets other individuals at a constant rate, the encounter resulting in an infection if the partner is susceptible. One may think that the rate of encounters in the local infections is not fixed, but depends on the number of members of the household. For example, translated in our notations, \cite{Neal2006} multiplies our rate of local infections by $\nu_i-1$, thus making the rate of encounters proportional to the size of the household. We expect that our results can be extended to that situation, probably at the cost of a stronger moment assumption on the random variables $\nu_i$.

This yields the following formal definition of the process.
Let
\begin{align*}
\X = \lbrace (n,k) \in \N \times \Z_+ : n \geq 1,\, 0 \leq k \leq n \rbrace.
\end{align*}

\begin{definition}[SIS household epidemic model] \label{def:household_model}
	Fix $ \lambda_L > 0 $, $ \lambda_G > 0 $ and $ \gamma > 0 $.
	Let $ \lbrace (\nu_i, X_i(0)), i \geq 1 \rbrace $ be i.i.d. $ \X $-valued random variables such that $ \E[\nu_1^2] < +\infty $ and let $ (P_{inf,i}(t), t \geq 0, i \geq 1) $ and $ (P_{rec,i}(t), t \geq 0, i \geq 0) $ be mutually independent standard Poisson processes, which are also independent of $ \lbrace (\nu_i, X_i(0)), i \geq 1 \rbrace $.
	We define $ \overline{\nu}^N = \frac{1}{N}\sum_{i=1}^{N} \nu_i $.
	For $ N \geq 1 $, let $ (X^N_1(t), \ldots, X^N_N(t), t \geq 0) $ be the solution of the following system of SDEs:
	\begin{multline}\label{Nsde}
	X_i^N(t)=X_i(0)-P_{rec,i}\left(\gamma\int_0^tX^N_i(s)ds\right)\\
	+P_{inf,i}\left(\int_0^t\left(1-\frac{X^N_i(s)}{\nu_i}\right)\left[\lambda_L \frac{\nu_i}{\nu_i-1} X^N_i(s)+\lambda_G \frac{\nu_i}{\overline{\nu}^N} \frac{1}{N} \sum_{j=1}^NX^N_j(s)\right]ds\right) ,
	\end{multline}
	for $ 1 \leq i \leq N $.
	We call this process the SIS household model with $ N $ households.
\end{definition}

The fact that there exists a unique solution to \eqref{Nsde} follows from a standard argument which exploits the fact that the jumps are isolated, and the process remains constant between its jumps.
The distribution of the $ \nu_i $'s will be fixed throughout the paper, and we set
\begin{align*}
\pi(n) = \P(\nu_1 = n), && \overline{\pi} = \E[\nu_1].
\end{align*}
We shall also use the size-biased distribution of the $ \nu_i $'s and its first moment, which we define as
\begin{align*}
\pi^+(n) = \frac{n\pi(n)}{\overline{\pi}}, && \overline{\pi}^+ = \sum_{n \geq 1} n \pi^+(n) = \frac{\E[\nu_1^2]}{\E[\nu_1]}.
\end{align*}

We note that the different households only interact through the mean number of infectious individuals in the $ N $ households, \textit{i.e.} it is a \textit{mean-field} interaction.
We thus expect that, as the number of households $ N $ becomes very large, any finite subset of households are asymptotically mutually independent and each one evolves according to the following SDE:
\begin{multline}\label{eq:MK}
X(t)=X(0)- P_{rec}\left(\gamma\int_0^tX(s)ds\right)\\
+P_{inf}\left(\int_0^t \left(1-\frac{X(s)}{\nu}\right) \left[\lambda_L \frac{\nu}{\nu-1} X(s)+\lambda_G \frac{\nu}{\overline{\pi}} \E[X(s)]\right]ds\right) ,
\end{multline}
where $(\nu, X(0))$ has the same law as $(\nu_1, X^N_1(0))$ and $ P_{inf} $ and $ P_{rec} $ are two independent standard Poisson processes which are also independent of $ (\nu, X(0)) $.
This is what is called \textit{propagation of chaos} \cite{ASS}, and will be made more precise in Theorem~\ref{thm:chaos} below.

This equation is a McKean-Vlasov Poisson driven SDE, because the transition rates of $ (X(t), t \geq 0) $ depend on the law of $ X(t) $ (specifically on its expectation).
We refer to McKean \cite{McK} for the study of similar Brownian driven SDEs.
As we will see later, this equation defines a semigroup acting on probability distributions on $ \N \times \Z_+ $ but, contrary to ordinary Markov processes, this semigroup is non-linear (because of the term $ \E[X(s)] $ appearing on the right hand side of \eqref{eq:MK}).
For this reason we will call $ (X(t), t \geq 0) $ the \textit{non-linear} Markov process.

\section{Main results} \label{sec:results}

\paragraph*{Existence and uniqueness of the non-linear process.}

It is not clear \textit{a priori} that there exists a process solving \eqref{eq:MK}, much less that it is unique.

Suppose for a moment that it exists and set
\begin{align*}
\mu_{n,k}(t) = \P(X(t) = k, \nu = n ).
\end{align*}
Then, $ \mu(t) = \lbrace \mu_{n,k}(t), (n,k) \in \X \rbrace $ is the law of $ (\nu, X(t)) $ and
\begin{align} \label{marginal_pi}
\forall n \geq 1, \quad \sum_{k=0}^{n} \mu_{n,k}(t) = \pi(n).
\end{align}
Equation~\eqref{eq:MK} then implies that $ \lbrace \mu_{n,k}(t), t \geq 0, (n,k) \in \X \rbrace  $ solves the following non-linear Fokker-Planck equation:
\begin{multline} \label{fokker_planck}
\frac{d\mu_{n,k}(t)}{dt}= \mu_{n,k-1}(t) \left(1-\frac{k-1}{n}\right) \left[ \lambda_L\frac{n}{n-1}(k-1)+\lambda_G\frac{n}{\overline{\pi}}\sum_{i=1}^\infty \sum_{j=1}^i j\mu_{i,j}(t) \right] \\ -\mu_{n,k}(t) \left\{ \left(1-\frac{k}{n}\right) \left[ \lambda_L\frac{n}{n-1} k+\lambda_G\frac{n}{\overline{\pi}}\sum_{i=1}^\infty \sum_{j=1}^i j\mu_{i,j}(t)\right] + \gamma k \right\} +\mu_{n,k+1}(t)\gamma (k+1),
\end{multline}
with the convention that $ \mu_{n,-1}(t) = \mu_{n,n+1}(t) = 0 $.
Note that \eqref{fokker_planck} defines an infinite system of coupled ordinary differential equations.
We then have the following theorem.

\begin{theorem}\label{thm:existence_uniqueness}
	Assume that the second moment of the probability distribution $\pi$ is finite.
	Then, given a probability measure $ \mu(0) = \{\mu_{n,k}(0), (n,k) \in \X \}$ satisfying \eqref{marginal_pi}, there exists a unique time dependent probability measure $(\mu(t),\ t\ge0 )$  on $ \X $ which solves the system of ODEs \eqref{fokker_planck}. 
	Moreover, given a random variable $ (\nu, X_0) $ which is such that $\P(X_0=k, \nu=n)=\mu_{n,k}(0)$ for $ (n,k) \in \X $, the SDE \eqref{eq:MK} has a unique solution $(X(t),\ t\ge0)$ which is such that for each $t\geq 0$, $\P(X(t)=k,\nu=n)=\mu_{n,k}(t)$ for each $ (n,k) \in \X $.
\end{theorem}

We prove this theorem in Section~\ref{sec:existence_uniqueness}.

\paragraph*{Propagation of chaos.}

We now deal with the limiting behaviour of the household model of Defintion~\ref{def:household_model} as the number of households $ N $ tends to infinity.
Let $ \mathcal{P}(D([0,\infty), \X)) $ denote the space of probability measures on the sample paths space $ D([0,\infty),\X) $, endowed with the topology of Skorokhod convergence on compact sets.
Also let $ \mu^\ast \in \mathcal{P}(D([0,\infty),\X)) $ denote the law of the non-linear Markov process $ ((\nu, X(t)), t \in [0,T]) $, given by Theorem~\ref{thm:existence_uniqueness}. We note that all elements $\mu$ of $ \mathcal{P}(D([0,\infty), \X)) $ which we shall consider are such that the first element $\nu$ of the pair $(\nu,X(t))$ is $\mu$ a.s. constant in time.

\begin{theorem}[Propagation of chaos in the SIS household model] \label{thm:chaos}
	Assume that $ \lbrace (\nu_i, X_i(0)), i \geq 1 \rbrace $ are independent and identically distributed $ \X $-valued random variables such that $ \E[\nu_1^2] < +\infty $.
	For all $ N \geq 1 $, let $ (X_i^N(t), t \geq 0, 1 \leq i \leq N) $ be the solution of equation~\eqref{Nsde}.
	Define $ \mu_N \in \mathcal{P}(D([0,\infty),\X)) $ by
	\begin{align*}
	\mu_N = \frac{1}{N} \sum_{i=1}^{N} \delta_{(\nu_i, X^N_i(\cdot))}.
	\end{align*}
	Then the random measure $ \mu_N $ converges weakly to $ \mu^* $ as $ N \to \infty $ in probability.
	Moreover, for any $ k \geq 1 $,
	\begin{align*}
	\text{Law}\left( (\nu_1, X_1^N(\cdot)), \ldots, (\nu_k, X_k^N(\cdot) \right) \Rightarrow (\mu^*)^{\otimes k} \: \text{ as } N \to \infty
	\end{align*}
	in $ \mathcal{P}\left(D([0,\infty),\X^k)\right) $.
\end{theorem}

We prove Theorem~\ref{thm:chaos} in Section~\ref{sec:chaos}.
Note that by Proposition~2.2 in \cite{ASS}, the second part of the theorem follows from the convergence of the empirical measures $ \mu_N $.
Theorem~\ref{thm:chaos} says two things: as $ N $ becomes large, any finite subset of households behaves asymptotically as independent copies of the non-linear Markov process \eqref{eq:MK}, and the global epidemic, as measured through the empirical measure $ \mu_N $, becomes asymptotically deterministic and equal to the law of the non-linear Markov process.
It is then natural to ask whether the epidemic has an endemic equilibrium and if it is stable in the non-linear Markov process.

In Theorem~\ref{thm:chaos}, the assumption that $ \{(\nu_i, X_i(0)), i \geq 1\} $ are \textit{i.i.d.} can be replaced by the weak convergence in distribution of the empirical measures
\begin{align*}
\frac{1}{N} \sum_{i=1}^{N} \delta_{(\nu_i, X^N_i(0))}
\end{align*}
to some deterministic measure $ \mu_0 $ (for example it is possible to assume that these empirical measures are deterministic and given by some sequence converging to $ \mu_0 $).
This does not affect the propagation of chaos result, and all our proofs carry over to that case.
Assuming that the initial condition is deterministic mainly reduces the fluctuations of $ \mu^N $ around $ \mu^* $.
This will be exposed in an upcoming paper treating the central limit theorem for the sequence $ (\mu^N, N \geq 1) $.

\paragraph*{Large time behaviour of the non-linear Markov process.}

As is usual in SIS epidemic models, there is in our model a basic reproduction number $ R_\ast $ such that if $ R_\ast > 1 $, there exists a unique stable endemic equilibrium (\textit{i.e.} the epidemic survives forever) and if $ R_\ast \leq 1 $, the disease free equilibrium is globally asymptotically stable (the epidemic eventually dies out).
This number is usually defined as the number of secondary infections produced by a single infectious individual.
Here, however, this number will be defined as the mean number of households which are infected by a single household, in which there is initially one infected individual and whose size is chosen according to the size-biased distribution $ \pi^+ $.

To do this, let $ (I(t), t \geq 0) $ be the solution to the following SDE:
\begin{align} \label{def:I}
I(t) = I(0) + P_{inf}\left( \int_{0}^{t} \lambda_L \left( 1 - \frac{I(s)-1}{\nu-1} \right) I(s) ds \right) - P_{rec}\left( \int_{0}^{t} \gamma I(s) ds \right),
\end{align}
where $ \nu $ is distributed according to the probability distribution $ \pi $ and $ P_{inf} $ and $ P_{rec} $ are two independent standard Poisson processes, which are independent of $ (\nu, I(0)) $.
Then $ (I(t), t \geq 0) $ is the number of infectious individuals in an isolated household.

We then define
\begin{align} \label{def_R0}
R_\ast = \lambda_G \sum_{n = 1}^{\infty} \pi^+(n) \E \left[ \left. \int_{0}^{\infty} I(t) dt \right| I(0) = 1, \nu = n \right].
\end{align}
The large time behaviour of the non-linear process $ (X(t), t \geq 0) $ of Theorem~\ref{thm:existence_uniqueness} is then given by the following result.

\begin{theorem}[Large time behaviour of the non-linear Markov process] \label{thm:large_time}
	Let $ (X(t), t \geq 0) $ be the unique solution to equation~\eqref{eq:MK}, and assume that $ \E[\nu^2] < + \infty $.
	\begin{longlist}
		\item If $ R_\ast > 1 $, then there exists a unique probability distribution $ \mu_\infty $ on $ \X $ such that, if $\, \P(X(0) \geq 1) > 0 $, $ (\nu, X(t)) $ converges in distribution to $ \mu_\infty $ as $ t \to \infty $.
		Moreover $ \mu_\infty $ is non-trivial in the sense that $ \mu_\infty \neq \pi \otimes \delta_0 $.
		\item If $ R_\ast \leq 1 $, then $ X(t) \to 0 $ in probability as $ t \to \infty $.
	\end{longlist}
\end{theorem}

We prove Theorem~\ref{thm:large_time} in Section~\ref{sec:large_time}.
This result should be seen as an analogue of the fact that the solution of the ODE
\begin{align} \label{ode_sis}
\frac{di(t)}{dt} = \lambda i(t)\left(1-\frac{i(t)}{n}\right) - \gamma i(t)
\end{align}
converges as $ t \to \infty $ to $ n\left(1-\frac{\gamma}{\lambda}\right) $ if $ \lambda > \gamma $ and to 0 otherwise.

In the proof of Theorem~\ref{thm:large_time}, we shall also prove the following formula for $ R_\ast $, which is of independent interest:
\begin{align*}
R_\ast = \frac{\lambda_G}{\gamma} \sum_{n=1}^\infty \pi^+(n) \left( 1 + \sum_{\ell=1}^{n-1} \left( \frac{\lambda_L}{\gamma} \right)^{\ell} \prod_{j=1}^\ell \left( 1 - \frac{j-1}{n-1} \right) \right),
\end{align*}
which can be rewritten as 
\begin{align} \label{formula_R0}
R_\ast=  \frac{\lambda_G}{\gamma} \sum_{n=1}^\infty \pi^+(n) \left(\frac{\lambda}{\gamma(n-1)}\right)^{n-1} (n-1)!
\sum_{k=0}^{n-1}\frac{1}{k!}\left(\frac{\gamma(n-1)}{\lambda_L}\right)^k
\end{align}
The first formula above appears in the proof of Lemma~\ref{lemma:R0} in Subsection~\ref{subsec:R0}. To deduce the second formula from the first, just note that $\prod_{j=1}^\ell\frac{n-j}{n-1}=(n-1)^{-\ell}\frac{(n-1)!}{(n-\ell-1)!}$ and let $k=n-\ell-1$. 
Note that similar formulas appear in \cite{Ball} and \cite{Neal2006}, with minor differences due to the fact that local infections follow slightly different rules in their models. 

\begin{remark}
	\begin{longlist}
		\item If $ \pi = \delta_1 $, every household is of size 1, and \eqref{Nsde} reduces to the homogeneous SIS epidemic model, with parameters $ \lambda_G $ and $ \gamma $ (see \cite{BS}).
		We can then check that $ \E[X(t)] $ solves the ODE \eqref{ode_sis} with $ \lambda = \lambda_G $, and that \eqref{def_R0} reduces to $ R_\ast = \lambda_G / \gamma $, as expected.
		\item The same is true if we take $ \lambda_L = 0 $ and keep $ \pi $ very general, the only infections in the system are global infections and the model reduces to the standard SIS epidemic model.
		\item Another interesting case is when the size of all the households is very large.
		In that case, if we approximate $ (I(t), t \geq 0) $ by a branching process, we see that $ R_\ast $ should be approximated by $ +\infty $ if $ \lambda_L \geq \gamma $ and by $ \lambda_G/(\gamma - \lambda_L) $ if $ \lambda_L < \gamma $,	and we see that $ R_\ast > 1 $ is equivalent to $ \lambda_G + \lambda_L > \gamma $.
	\end{longlist}
\end{remark}

\section{The forced process} \label{sec:forced_process}

It is worth noting that if we replace $ \E[X(s)] $ in \eqref{eq:MK} by any deterministic measurable function $ s \mapsto m(s) $, then $ (X(t), t \geq 0) $ becomes a time-inhomogeneous Markov process.

\begin{definition}[The forced process]
	Let $ m : \R_+ \to [0,\bar{\pi}] $ be a measurable function, $ (\nu, X_0) $ an $ \X $-valued random variable and $ P_{inf} $ and $ P_{rec} $ two independent standard Poisson processes which are also independent of $ (\nu, X_0) $.
	Then the forced process $ (X_t(m), t \geq 0) $ is defined as the solution to
	\begin{multline} \label{eq:X_m}
	X_t(m) = X_0 + P_{inf}\left( \int_{0}^{t} \left[ \lambda_L\frac{\nu}{\nu-1} X_s(m)  + \lambda_G \frac{\nu}{\pib}m(s) \right]\left(1-\frac{X_s(m)}{\nu}\right) ds \right) \\ - P_{rec}\left( \int_{0}^{t} \gamma X_s(m) ds \right).
	\end{multline}
\end{definition}

We call this process the \textit{forced} process because we fix the intensity of global infections to be $ \lambda_G\nu m(t)/\pib $.
The fact that there exists a unique strong solution to \eqref{eq:X_m} follows from standard arguments similar to that used in \eqref{Nsde}.

Comparing \eqref{eq:MK} and \eqref{eq:X_m}, we see that solving \eqref{eq:MK} is equivalent to finding a measurable function $ m $ such that $ m(t) = \E[X_t(m)] $ for all $ t \geq 0 $.

\subsection{Graphical construction of the forced process}

We are going to show that we can construct this process with the following procedure.
The same construction, which was first introduced by Harris \cite{Harris} in a different context,  was used in \cite{Neal2008} for an SIS model on a large circle with local and global infectious contacts.
Let $ c(dk) $ denote the counting measure on $ \N $.
Conditionally on $ (\nu, X_0) $, let $ \Pi_{rec} $, $ \Pi_{L} $ and $ \Pi_{G} $ be three mutually independent Poisson point processes such that
\begin{itemize}
	\item $ \Pi_{rec} $ is a Poisson point process on $ \R_+ \times \llbracket 1, \nu \rrbracket $ with intensity $ \gamma dt \otimes c(dk) $,
	\item $ \Pi_L $ is a Poisson point process on $ \R_+ \times \llbracket 1, \nu \rrbracket \times \llbracket 1, \nu \rrbracket $ with intensity $ \frac{\lambda_L}{\nu-1} dt \otimes c(dk) \otimes c(dk) $,
	\item $ \Pi_G $ is a Poisson point process on $ \R_+ \times \llbracket 1, \nu \rrbracket \times [0,\overline{\pi}] $ with intensity $ \frac{\lambda_G}{\pib} dt \otimes c(dk) \otimes du $.
\end{itemize}
Let us describe the effect of these different processes before formally constructing the forced process.
A point $ (t, i) $ in $ \Pi_{rec} $ means that if the individual $ i $ is infected at time $ t^- $, it becomes susceptible at time $ t $ (it undergoes a remission).
A point $ (t, i, j) $ in $ \Pi_{L} $ means that if individual $ i $ is infectious and individual $ j $ is susceptible at time $ t^- $, then individual $ j $ becomes infected at time $ t $.
Finally a point $ (t,i,u) $ in $ \Pi_G $ means that a global infection targetting individual $ i $ takes place at time $ t $ if individual $ i $ is susceptible at time $ t^- $ and if $ u \leq m(t) $.

In fact, we can view the total set of infectious individuals at any time as the union of several local epidemics, each resulting from a previous global infection or from the individuals infected at time 0.
To do this, note that $ \Pi_G $ is almost surely locally finite, so we can order its points according to their time coordinate.
Thus let
\begin{align*}
\Pi_G = \lbrace (t_k, i_k, u_k), k \geq 1, 0 < t_1 < t_2 < \ldots \rbrace.
\end{align*}
Let us then define a random set $ I^k(t) \subset \llbracket 1, \nu \rrbracket $ for all $ t \geq 0 $ as follows.
\begin{itemize}
	\item For $ t < t_k $, $ I^k(t) = \emptyset $.
	\item At $ t = t_k $, we set $ I^k(t_k) = \lbrace i_k \rbrace $.
	\item For each $ (t,i,j) \in \Pi_L $, if $ i \in I^k(t^-) $, then $ I^k(t) = I^k(t^-) \cup \lbrace j \rbrace $.
	\item For each $ (t, i) \in \Pi_{rec} $, $ I^k(t) = I^k(t^-) \cap \lbrace i \rbrace^c $.
\end{itemize}
We define in the same way the local epidemic resulting from the initially infected individuals $ (I^0(t), t \geq 0) $, \textit{i.e.} $ I^0(0) = \lbrace i : 1 \leq i \leq X_0 \rbrace $ and $ I^0 $ evolves according to the same rules as $ I^k $ for $ k \geq 1 $.
We note that, for all $ k \geq 0 $, $ (I^k(t), t \geq 0) $ is right-continuous with left limits.

\begin{proposition} \label{prop:forced_process}
	For all $ t \geq 0 $, let
	\begin{align} \label{def_X_m}
	X_t(m) = \left| I^0(t) \cup \bigcup_{k \geq 1} \lbrace I^k(t) : u_k \leq m(t_k) \rbrace \right|,
	\end{align}
	where $ |\cdot| $ denotes the cardinal of a set.
	Then the process $ (X_t(m), t \geq 0) $ is a solution to the SDE \eqref{eq:X_m}.
\end{proposition}

\begin{proof}
	Clearly $ X_0(m) = | I^0(0) | = X_0 $.
	It remains to check that the waiting times between upward and downward jumps of $ X_t(m) $ are distributed as exponential variables with the correct rates.
	
	If the current value of $ X_t(m) $ is $ x $, then the next remission takes place at the time
	\[ \inf\{ s>t:\ \exists i \in \llbracket 1,\nu\rrbracket, k\ge 0 \text{ s.t. } (s,i)\in \Pi_{rec},\ i\in I^k(s^-) \text{ and }u_k \leq m(t_k) \}\,.\]
	(Note that we can set $ u_0 = \overline{\pi} $ and $ t_0 = 0 $.)
	This happens at instantaneous rate $ \gamma x $, as in \eqref{eq:X_m}.
	
	Likewise, the first time after $t$ that a susceptible individual becomes infected due to a local infcetion is given by the next point $ (s, i, j) \in \Pi_{L} $ with $ s > t $ such that $ i \in I^k(s^-) $ for some $ k \geq 0 $ with $ u_k \leq m(t_k) $ and $ j \notin I^k(s^-) $.
	This happens at rate $ \frac{\lambda_L}{\nu-1} x (\nu - x) $, as in \eqref{eq:X_m}.
	
	Finally, the next time a susceptible individual becomes infected due to a global infection is the next $ (s, i, u) \in \Pi_G $ with $ s > t $ such that $ i \notin I^k(s^-) $ for all $ k \geq 0 $ such that $ u_k \leq m(t_k) $ and $ u \leq m(s) $.
	This happens at instantaneous rate $ \frac{\lambda_G}{\pib} (\nu - x) m(t) $, as in \eqref{eq:X_m}.
\end{proof}

\subsection{Monotonicity of the forced process}

With this construction, the next lemma is straightforward.

\begin{lemma}[Monotonicity of the forced process] \label{lemma:monotonicity}
	Suppose that $ X_0^{(1)} $ and $ X_0^{(2)} $ are two random variables such that $ X_0^{(1)} \leq X_0^{(2)} $ almost surely.
	Also let $ m_1 $ and $ m_2 $ be two measurable functions from $ \R_+ $ to $ [0,\overline{\pi}] $ such that $ m_1(t) \leq m_2(t) $ for almost every $ t \geq 0 $.
	Then there exists a process $ (X_t(m_1), t \geq 0) $ solving \eqref{eq:X_m} with $ m = m_1 $ and $ X_0 = X_0^{(1)} $, and a process $ (X_t(m_2), t \geq 0) $ solving \eqref{eq:X_m} with $ m = m_2 $ and $ X_0 = X_0^{(2)} $, defined on the same probability space, such that, almost surely,
	\begin{align*}
	X_t(m_1) \leq X_t(m_2), \quad \forall t \geq 0.
	\end{align*}
\end{lemma}

\begin{proof}
	We use Proposition~\ref{prop:forced_process} to construct both processes with the same Poisson point processes $ \Pi_{rec} $, $ \Pi_L $ and $ \Pi_G $.
	We define $ (I^{0,i}(t), t \geq 0) $ for $ i \in \lbrace 1, 2 \rbrace $ as above with
	\begin{align*}
	I^{0,i}(0) = \lbrace k : 1 \leq k \leq X_0^{(i)} \rbrace,
	\end{align*}
	so that, almost surely, $ I^{0,1}(0) \subset I^{0,2}(0) $.
	Then, from the evolution of $ (I^{0,i}(t), t \geq 0) $, we deduce that $ I^{0,1}(t) \subset I^{0,2}(t) $ for all $ t \geq 0 $.
	Furthermore, since $ m_1 \leq m_2 $,
	\begin{align*}
	\lbrace k : u_k \leq m_1(t_k) \rbrace \subset \lbrace k : u_k \leq m_2(t_k) \rbrace.
	\end{align*}
	It then follows from equation~\eqref{def_X_m} that $ X_t(m_1) \leq X_t(m_2) $, where $ X_t(m_1)$ and $X_t(m_2) $ are constructed as in Proposition \ref{prop:forced_process}.
\end{proof}

The following lemma will also be useful in the proof of existence and uniqueness of the non-linear process.
For $ t \geq 0 $, $ m : \R_+ \to [0,\overline{\pi}] $ measurable and $ \mu_0 $ a probability measure on $ \X $ whose first marginal is $ \pi $, let
\begin{align} \label{eq:def_mu_bar_m}
\overline{\mu}_t(m, \mu_0) = \E [X_t(m)],
\end{align}
where $ (\nu, X_0) $ is distributed according to $ \mu_0 $.

\begin{lemma} \label{lemma:bound_mu_bar}
	Suppose that $ \mu_0 $ is as above.
	If $ m_1 $ and $ m_2 $ are two measurable functions from $ \R_+ $ to $ [0,\overline{\pi}] $ satisfying $ m_1(t) \leq m_2(t) $ for almost every $ t \geq 0 $, then
	\begin{align*}
	0 \leq \overline{\mu}_t(m_2, \mu_0) - \overline{\mu}_t(m_1, \mu_0) \leq \pibp\lambda_G \int_{0}^{t} (m_2(s) - m_1(s)) ds.
	\end{align*}
\end{lemma}

\begin{proof}
	The fact that $ \overline{\mu}_t(m_2) - \overline{\mu}_t(m_1) \geq 0 $ follows from Lemma~\ref{lemma:monotonicity}.
	To prove the second inequality, we construct $ (X_t(m_1), t \geq 0) $ and $ (X_t(m_2), t \geq 0) $ as in Proposition~\ref{prop:forced_process}.
	Then
	\begin{align*}
	0 \leq X_t(m_2) - X_t(m_1) \leq \left| \cup_{k \geq 1} \lbrace I^k(t) : m_1(t_k) \leq u_k \leq m_2(t_k) \rbrace \right|.
	\end{align*}
	Moreover, we can restrict the union to the values of $ k $ for which $ t_k \leq t $.
	Since $ |I^k(t)| \leq \nu $ for all $ t \geq 0 $, we can write
	\begin{align} \label{eq:inequality_X_m}
	0 \leq X_t(m_2) - X_t(m_1) \leq \nu \left| \lbrace k \geq 1 : m_1(t_k) < u_k \leq m_2(t_k), t_k \leq t \rbrace \right|.
	\end{align}
	Now, by the definition of $ \Pi_G $, the right hand side is, conditionally on $ \nu $, $ \nu $ times a Poisson random variable with parameter
	\begin{align*}
	\lambda_G\frac{\nu}{\pib} \int_{0}^{t} (m_2(s) - m_1(s)) ds.
	\end{align*}
	As a result, taking expectations in \eqref{eq:inequality_X_m} (first conditionally on $ \nu $ and then over the law of $ \nu $), we obtain
	\begin{align*}
	0 \leq \overline{\mu}_t(m_2, \mu_0) - \overline{\mu}_t(m_1, \mu_0) \leq  \lambda_G\,\pibp \int_{0}^{t} (m_2(s) - m_1(s)) ds,
	\end{align*}
	and the lemma is proved.
\end{proof}

We shall come back to the forced process in the proof of Theorem~\ref{thm:large_time}, as it will be used to characterize the possible stationary distributions of the non-linear process.

\section{Existence and uniqueness of the non-linear Markov process} \label{sec:existence_uniqueness}

We now set out to prove Theorem~\ref{thm:existence_uniqueness}.
We note that finding a solution to \eqref{fokker_planck} is equivalent to finding a fixed point of
\begin{align} \label{fixed_point}
m(\cdot) \mapsto \overline{\mu}_\cdot(m,\mu_0).
\end{align}
Indeed, if $ m_* $ is a fixed point of this function, then $ (X_t(m_*), t \geq 0) $ is a solution to \eqref{eq:MK}.
We thus need to prove that, given $ \mu_0 $, there exists a unique fixed point of \eqref{fixed_point}.

\begin{proof}[Proof of Theorem~\ref{thm:existence_uniqueness}]
	Fix $ \mu_0 $ and assume that $ (\nu, X_0) $ is distributed according to $ \mu_0 $.
	Let $ (m^{+,k}, k \geq 0) $ and $ (m^{-,k}, k \geq 0) $ be two sequences of functions defined by
	\begin{align*}
	&m^{+,0}(t) = \overline{\pi}, \qquad m^{+,k+1}(t) = \overline{\mu}_t(m^{+,k}, \mu_0), \\
	&m^{-,0}(t) = 0, \qquad m^{-,k+1}(t) = \overline{\mu}_t(m^{-,k}, \mu_0),
	\end{align*}
	where $ \overline{\mu}_t(m,\mu_0) $ was defined in \eqref{eq:def_mu_bar_m}.
	Clearly, since $ 0 \leq \E[X_t(m)] \leq \overline{\pi} $,
	\begin{align*}
	m^{+,1}(t) \leq m^{+,0}(t), && m^{-,1}(t) \geq m^{-,0}(t).
	\end{align*}
	Then by induction, using Lemma~\ref{lemma:monotonicity}, we obtain that, for all $ k \geq 0 $,
	\begin{align*}
	m^{-,k}(t) \leq m^{-,k+1}(t) \leq m^{+,k+1}(t) \leq m^{+,k}(t).
	\end{align*}
	Hence $ m^{+,k} $ and $ m^{-,k} $ both converge pointwise.
	Let $ m^{+,\infty} $ and $ m^{-,\infty} $ be their respective limits.
	Then, using Lemma~\ref{lemma:bound_mu_bar} with $ m_1 = m^{+,\infty} $ and $ m_2 = m^{+,k} $,
	\begin{align*}
	\left| \overline{\mu}_t(m^{+,\infty}) - m^{+,\infty}(t) \right| \leq \left| m^{+,k+1}(t) - m^{+,\infty}(t) \right| +\lambda_G\,\pibp \int_{0}^{t} (m^{+,k}(s) - m^{+,\infty}(s)) ds.
	\end{align*}
	The integral on the right hand side vanishes as $ k \to \infty $ by dominated convergence and the first term vanishes because $ m^{+,k} $ converges pointwise to $ m^{+,\infty} $.
	As a result, $ m^{+,\infty} $ (and also $ m^{-,\infty} $ by the same argument) is a fixed point of \eqref{fixed_point}.
	This shows existence of solutions to \eqref{fokker_planck} (and thus to \eqref{eq:MK}).
	
	To prove uniqueness, first note that, by induction and using Lemma~\ref{lemma:monotonicity}, any fixed point $ m_* $ satisfies
	\begin{align*}
	m^{-,k}(t) \leq m_*(t) \leq m^{+,k}(t),
	\end{align*}
	for all $ k \geq 0 $ and $ t \geq 0 $.
	Hence we also have
	\begin{align*}
	m^{-,\infty}(t) \leq m_*(t) \leq m^{+,\infty}(t).
	\end{align*}
	To prove uniqueness, it is thus enough to prove that $ m^{+,\infty}(t) = m^{-,\infty}(t) $ for all $ t \geq 0 $.
	Using Lemma~\ref{lemma:bound_mu_bar} with $ m_1 = m^{-,0} $ and $ m_2 = m^{+,0} $, we obtain
	\begin{align*}
	0 \leq m^{+,1}(t) - m^{-,1}(t) \leq \pib\,\pibp \lambda_G t,
	\end{align*}
	and by induction, we deduce that, for $ k \geq 1 $,
	\begin{align*}
	0 \leq m^{+,k}(t) - m^{-,k}(t) \leq \overline{\pi} \frac{(\pibp \lambda_G t)^k}{k!}.
	\end{align*}
	Leting $ k \to \infty $, it follows that $ m^{+,\infty}(t) = m^{-,\infty}(t) $ for all $ t \geq 0 $ and the theorem is proved.
\end{proof}

\section{Propogation of chaos for the SIS household model} \label{sec:chaos}

The aim of this section is to prove Theorem~\ref{thm:chaos}.
As we have explained before, using Proposition~2.2 in \cite{ASS}, the second part of the statement follows from the convergence of the empirical measures $ \mu_N $ to the law of the non-linear process $ \mu^* $.
We establish this convergence by showing that the sequence $ \lbrace \mu_N, N \geq 1 \rbrace $ is tight in $ \mathcal{P}(D([0,\infty),\X)) $, and identifying its possible limit points.

\begin{lemma} \label{lemma:tightness}
	The sequence $ \lbrace \mu_N, N \geq 1 \rbrace $ is tight in $ \mathcal{P}(D([0,\infty),\X)) $.
\end{lemma}

\begin{proof}
	By Proposition~2.2(ii) in \cite{ASS}, the sequence $ \lbrace \mu_N, N \geq 1 \rbrace $ is tight if and only if the laws of $ (\nu_1, X_1^N(\cdot)) $ are tight, but this is straightforward from \eqref{Nsde} where we see that the rate of increase is bounded by 
	$(\lambda_L+\lambda_G)\nu_i$, while the rate of decrease is bounded by $\gamma\nu_i$.
	%
\end{proof}

Next we note that equation~\eqref{Nsde} can be reformulated as follows.
Let $ \lbrace \M_{inf,i}, i \geq 1 \rbrace$ and $ \lbrace \M_{rec,i}, i \geq 1 \rbrace $ be mutually independent random Poisson measures on $\R_+^2$ with intensity measure the Lebesgue measure, which are also independent of $ \lbrace (\nu_i,X_i(0)), i \geq 1 \rbrace $.
Then, with the notations
\begin{align*}
\overline{\mu}^N_t := \frac{1}{N} \sum_{i=1}^{N} X_i^N(t),\quad
\Lambda(n, x, y, m) =: \left( 1 - \frac{x}{n} \right) \left[ \lambda_L \frac{n}{n-1} x + \lambda_G \frac{n}{y} m \right], 
\end{align*}
\begin{align*}
X_i^N(t)&= X_i(0) + \int_0^t \int_0^\infty {\bf1}_{u \le \Lambda(\nu_i, X^N_i(s^-), \overline{\nu}^N, \overline{\mu}^N_{s^-}) } \M_{inf,i}(ds\, du)\\
&\quad - \int_0^t \int_0^\infty {\bf1}_{u \le \gamma X^N_i(s^-)} \M_{rec,i}(ds\, du).
\end{align*}
Clearly, for any $\phi:\X\mapsto\R$,
\begin{multline*}
\phi(\nu_i,X^N_i(t)) = \phi(\nu_i,X^N_i(0))\\
+ \int_0^t \int_0^\infty \left[\phi(\nu_i, X^N_i(s^-) + 1) - \phi(\nu_i, X^N_i(s^-))\right] {\bf1}_{u \le \Lambda(\nu_i, X^N_i(s^-), \overline{\nu}^N, \overline{\mu}^N_{s^-})} \M_{inf,i}(ds,du)\\
+ \int_0^t \int_0^\infty \left[\phi(\nu_i, X^N_i(s^-)-1) - \phi(\nu_i, X^N_i(s^-))\right] {\bf1}_{u \le \gamma X^N_i(s^-)} \M_{rec,i}(ds,du).
\end{multline*}
Let $\MM_{inf,i}$ and $\MM_{rec,i}$ denote the compensated measures
\begin{align*}
&\MM_{inf,i}(ds,du)=\M_{inf,i}(ds,du)-dsdu, \\
&\MM_{rec,i}(ds,du)=\M_{rec,i}(ds,du)-dsdu.
\end{align*}
Then, setting
\begin{multline*}
M^\phi_i(t) = \int_0^t \int_0^\infty \left[ \phi(\nu_i, X^N_i(s^-)-1) - \phi(\nu_i, X^N_i(s^-)) \right] {\bf1}_{u \le \gamma X^N_i(s^-)} \MM_{rec,i}(ds,du) \\ + \int_0^t \int_0^\infty \left[ \phi(\nu_i, X^N_i(s^-)+1) - \phi(\nu_i, X^N_i(s^-)) \right] {\bf1}_{u \le \Lambda(\nu_i, X^N_i(s^-), \overline{\nu}^N, \overline{\mu}^N_{s^-})} \MM_{inf,i}(ds,du),
\end{multline*}
we have
\begin{multline*}
\phi(\nu_i,X^N_i(t)) = \phi(\nu_i,X^N_i(0)) + \gamma \int_0^t \left[ \phi(\nu_i, X^N_i(s)-1) - \phi(\nu_i, X^N_i(s)) \right] X^N_i(s) ds\\ + \int_0^t \left[ \phi(\nu_i, X^N_i(s)+1) - \phi(\nu_i, X^N_i(s)) \right] \Lambda(\nu_i, X^N_i(s), \overline{\nu}^N, \overline{\mu}^N_{s}) ds + M^\phi_i(t).
\end{multline*}
We rewrite this identity in the form
\begin{equation}\label{eq:phi_i}
\phi(\nu_i,X^N_i(t))=\phi(\nu_i,X^N_i(0))+\int_0^t[\mathcal L\phi](\nu_i,X^N_i(s),\overline{\nu}^N,\overline{\mu}^N_s)ds+M^\phi_i(t),
\end{equation}
where for $n\ge1$, $x\in\{0,1,\ldots,n\}$, $y\ge0$ and $ 0 \leq m \leq y $,
\begin{align*}
\mathcal L\phi(n,x,y,m) &= [\phi(n,x+1)-\phi(n,x)] \left(1-\frac{x}{n}\right) \left[\lambda_L\frac{n}{n-1}x+\lambda_G\frac{n}{y}m\right]\\&\quad + [\phi(n,x-1)-\phi(n,x)] \gamma x.
\end{align*}

\begin{proof}[Proof of Theorem~\ref{thm:chaos}]
	Let $ \mu^\infty $ be a limit point of the sequence $ \mu^N $.
	First note that, by the classical law of large numbers, for any bounded and measurable $ \phi : \X \to \R $,
	\begin{align*}
	\E_{\mu^\infty}[\phi(\nu,X(0))] = \E[\phi(\nu_1,X_1(0))].
	\end{align*}
	In order to identify the possible limit points of $ \mu^N $, we define, for $ \mu \in \mathcal{P}(D([0,\infty),\X)) $ and $ 0 \leq s \leq t \leq T $,
	\begin{align*}
	\Phi_{s,t}(\mu) = \E_\mu \left[ \left( \phi(\nu, X(t)) - \phi(\nu, X(s)) - \int_{s}^{t} \mathcal{L}\phi(\nu, X(r), \overline{\nu}, \overline{\mu}_r) dr \right) \psi_s(\nu, X(\cdot)) \right],
	\end{align*}
	where
	\begin{align*}
	\overline{\mu}_t = \E_\mu \left[ X(t) \right], && \overline{\nu} = \E_\mu[\nu],
	\end{align*}
	and $ \phi $ is any bounded function from $ \X $ to $ \R $ and $ \psi_s $ is of the form
	\begin{align*}
	\psi_s(\nu, X(\cdot)) = \phi_1(\nu, X(s_1))\ldots\phi_k(\nu,X(s_k))
	\end{align*}
	with $ 0 \leq s_1 \leq \ldots \leq s_k \leq s $ and $ \phi_1, \ldots, \phi_k $ are bounded functions from $ \X $ to $ \R $.
	By Theorem~\ref{thm:existence_uniqueness}, the result will be proved if we show that
	\begin{align*}
	\Phi_{s,t}(\mu^\infty) = 0, 
	\end{align*}
	almost surely for any such function $ \Phi_{s,t} $.
	Using \eqref{eq:phi_i},
	\begin{align*}
	\Phi_{s,t}(\mu^N) = \frac{1}{N} \sum_{i=1}^{N} (M^\phi_i(t) - M^\phi_i(s)) \psi_s(\nu_i, X_i^N(\cdot)).
	\end{align*}
	From the definition of $ M^\phi_i $,
	\begin{align*}
	\langle M^\phi_i, M^\phi_j \rangle_t = 0, \quad \forall i \neq j,
	\end{align*}
	and
	\begin{align*}
	\langle M^{\phi}_i \rangle_t = \int_{0}^{t} \mathcal{G}\phi(\nu_i, X^N_i(s),\overline{\nu}^N, \overline{\mu}^N_s) ds,
	\end{align*}
	where
	\begin{align*}
	\mathcal{G}\phi(n,x,y,m) &= \left[ \phi(n,x+1)-\phi(n,x) \right]^2 \left( 1 - \frac{x}{n} \right) \left[\lambda_L\frac{n}{n-1} x + \lambda_G \frac{n}{y}m\right] \\ & \quad+ \left[ \phi(n,x-1) - \phi(n,x) \right]^2 \gamma x.
	\end{align*}
	Note that, for $ m \leq y $
	\begin{align*}
	\mathcal{G}\phi(n,x,y,m) \leq 4 \sup_\X |\phi|^2 \left( \lambda_L + \lambda_G + \gamma \right) n.
	\end{align*}
	As a result,
	\begin{align*}
	\E\left[ \Phi_{s,t}(\mu^N)^2 \right] &= \frac{1}{N^2} \sum_{i=1}^{N} \E \left[ (\langle M^\phi_i \rangle_t - \langle M^\phi_i \rangle_s) \psi_s(\nu_i, X^N_i(\cdot))^2 \right] \\
	& \leq \frac{C}{N} \E[\nu_1],
	\end{align*}
	for some $ C > 0 $.
	It follows that
	\begin{align*}
	\Phi_{s,t}(\mu^N) \to 0,
	\end{align*}
	in $ L^2 $ as $ N \to \infty $, hence $ \mu^\infty $ is equal to $ \mu^\ast $, the distribution of the non-linear process of \eqref{eq:MK}.
	This proves Theorem~\ref{thm:chaos}.
\end{proof}

\section{Large time behaviour of the non-linear Markov process} \label{sec:large_time}

Let us start this section by noting that if the non-linear process of \eqref{eq:MK} with initial distribution $ \mu_0 $ is stationary, then the forced process with initial distribution $ \mu_0 $ and with $ m(t) = \overline{\mu}_0 = \E_{\mu_0}[X(0)] $ is also stationary.
Thus to study the possible stationary distributions of the non-linear process, we first study the large-time behaviour of the forced process.

\subsection{The large-time behaviour of the forced process}

Suppose that we take $ m(t) = m $ for all $ t \geq 0 $ for some $ m \in [0, \overline{\pi}] $.
Then $ (X_t(m), t \geq 0) $ becomes a homogeneous continuous-time Markov process.
On the event $ \lbrace \nu = n \rbrace $, it takes values in $ \llbracket 0, n \rrbracket $.
If $ m > 0 $, then it is positive recurrent on this set, while if $ m = 0 $, 0 is the only absorbing state for $ X_t(m) $.
As a result, conditionally on $ \nu = n $, $ (X_t(m), t \geq 0) $ admits a unique stationary distribution.
It follows that $ ((\nu, X_t(m)), t \geq 0) $ admits a unique stationary distribution $ \mu_\infty(m) $.

This distribution can be obtained as in Proposition~\ref{prop:forced_process} in the following way.
Let $ \overleftarrow{\Pi}_{rec} $, $ \overleftarrow{\Pi}_L $ and $ \overleftarrow{\Pi}_G $ be independent Poisson point processes as above, but on $ \R_- $ instead of $ \R_+ $ for the first coordinate.
We can then order the points in $ \overleftarrow{\Pi}_G $ in decreasing order:
\begin{align*}
\overleftarrow{\Pi}_G = \lbrace (t_\ell, i_\ell, u_\ell), \ell \geq 1, 0 > t_1 > t_2 > \ldots \rbrace.
\end{align*}
The points in $ \overleftarrow{\Pi}_G $ represent global infection which took place in the past.
We then perform the same construction of $ I^k(t) $, this time for $ t \leq 0 $, and we set
\begin{align*}
X_\infty(m) = \left| \cup_{k \geq 1} \lbrace I^k(0) : u_k \leq m \rbrace \right|.
\end{align*}

\begin{proposition} \label{prop:X_inf}
	For each $ m \in [0,\overline{\pi}] $, $ (\nu, X_\infty(m)) $ is distributed according to $ \mu_\infty(m) $.
\end{proposition}

\begin{proof}
	For $ t \geq 0 $, let
	\begin{align*}
	\tilde{X}_t(m) = \left| \cup_{k \geq 1} \lbrace I^k(0) : u_k \leq m, t_k \geq -t \rbrace \right|.
	\end{align*}
	In other words, we only consider the local epidemics which started after time $ -t $.
	Then from Proposition~\ref{prop:forced_process}, we see that for each $ t \geq 0 $, $ \tilde{X}_t(m) $ is distributed as $ X_t(m) $, where $ X_t(m) $ is the solution of \eqref{eq:X_m} with $ m(t) = m $ and $ X_0 = 0 $.
	By the ergodic theorem for homogeneous Markov processes, $ (\nu, X_t(m)) $, and hence $ (\nu, \tilde{X}_t(m)) $, converge in distribution as $ t \to \infty $ to $ \mu_\infty(m) $.
	At the same time, we see from the definition of $ \tilde{X}_t(m) $ and $ X_\infty(m) $ that
	\begin{align*}
	\tilde{X}_t(m) = \sum_{i=1}^{\nu} \indic{\exists k \geq 1 : i \in I^k(0), u_k \leq m, t_k \geq -t}, && X_\infty(m) = \sum_{i=1}^{\nu} \indic{\exists k \geq 1 : i \in I^k(0), u_k \leq m}.
	\end{align*}
	Hence by monotone convergence,
	\begin{align*}
	\tilde{X}_t(m) \rightarrow X_\infty(m) \text{ as } t \to \infty,
	\end{align*}
	almost surely, and the lemma is proved.
\end{proof}

The next lemma says that $ m(t) \to m_\infty $ as $ t \to \infty $ is sufficient for $ X_t(m) $ to converge in distribution to $ \mu_\infty(m_\infty) $.

\begin{lemma}[Large time behaviour of the forced process] \label{lemma:limit_forced_process}
	Suppose that $ m : \R_+ \to [0,\overline{\pi}] $ is measurable and that
	\begin{align*}
	m_\infty = \lim_{t \to \infty} m(t)
	\end{align*}
	exists.
	Then $ (\nu, X_t(m)) $ converges in distribution to $ \mu_\infty(m_\infty) $ as $ t $ tends to infinity.
\end{lemma}

\begin{proof}
	Suppose for now that $ 0 < m_\infty < \overline{\pi} $.
	Then for all $ \varepsilon > 0 $, there exists $ t_\varepsilon $ such that, for all $ t \geq t_\varepsilon $,
	\begin{align*}
	m_\infty - \varepsilon \leq m(t) \leq m_\infty + \varepsilon.
	\end{align*}
	We choose $ \varepsilon $ small enough that $ 0 \leq m_\infty - \varepsilon $ and $ m_\infty + \varepsilon \leq \overline{\pi} $.
	We then define two functions $ m^+ $ and $ m^- $ by
	\begin{align*}
	m^+(t) = \overline{\pi}\indic{t < t_\varepsilon} + (m_\infty + \varepsilon)\indic{t \geq t_\varepsilon}, && m^-(t) = (m_\infty - \varepsilon)\indic{t \geq t_\varepsilon}.
	\end{align*}
	Then $ m^- \leq m \leq m^+ $, so by Lemma~\ref{lemma:monotonicity}, we can construct jointly the three processes $ (X_t(m^-), t \geq 0) $, $ (X_t(m), t \geq 0) $ and $ (X_t(m^+), t \geq 0) $ such that, almost surely,
	\begin{align*}
	X_t(m^-) \leq X_t(m) \leq X_t(m^+), \quad \forall t \geq 0.
	\end{align*}
	It follows that, for each $ t \geq 0 $ and each $ (n,k) \in \X $,
	\begin{align*}
	\P\left( \nu = n, X_t(m^+) \leq k \right) \leq \P\left( \nu = n, X_t(m) \leq k \right) \leq \P\left( \nu = n, X_t(m^-) \leq k \right).
	\end{align*}
	Since $ m^+ $ and $ m^- $ are both constant after time $ t_\varepsilon $ (which is deterministic), as $ t \to \infty $, $ (\nu, X_t(m^+)) $ and $ (\nu, X_t(m^-)) $ respectively converge in distribution to $ \mu_\infty(m_\infty + \varepsilon) $ and $ \mu_\infty(m_\infty - \varepsilon) $.
	Thus, letting $ t \to \infty $ above,
	\begin{multline*}
	\mu_\infty(m_\infty + \varepsilon)\left( \lbrace n \rbrace \times \lbrace 0, \ldots, k \rbrace \right) \leq \liminf_{t \to \infty} \P\left( \nu = n, X_t(m) \leq k \right) \\ \leq \limsup_{t \to \infty} \P \left( \nu = n, X_t(m) \leq k \right) \leq \mu_\infty(m_\infty - \varepsilon)\left( \lbrace n \rbrace \times \lbrace 0, \ldots, k \rbrace \right).
	\end{multline*}
	But, as $ \varepsilon \downarrow 0 $, the measure $ \mu_\infty(m_\infty \pm \varepsilon) $ converges weakly to $ \mu_\infty(m_\infty) $ (in fact the construction in Proposition~\ref{prop:X_inf} gives a construction of $ X_\infty(m \pm \varepsilon) $ and $ X_\infty(m) $ such that $ X_\infty(m\pm \varepsilon) \to X_\infty(m) $ almost surely as $ \varepsilon \downarrow 0 $, using monotone convergence as in the proof of Proposition~\ref{prop:X_inf}).
	Hence letting $ \varepsilon \downarrow 0 $ above, we obtain, for any $ (n,k) \in \X $,
	\begin{align*}
	\P \left( \nu = n, X_t(m) \leq k \right) \rightarrow \mu_\infty(m_\infty)\left( \lbrace n \rbrace \times \lbrace 0, \ldots, k \rbrace \right) \text{ as } t \to \infty,
	\end{align*}
	and the lemma is proved.
	If $ m_\infty = 0 $, then we can take instead $ m^-(t) = 0 $, and if $ m_\infty = \overline{\pi} $, then we take $ m^+(t) = \overline{\pi} $, and the rest of the proof is essentially identical.
\end{proof}

\subsection{The stationary distribution of the forced process}

We now study in more detail the family of distributions $ \mu_\infty(\cdot) $.
For $ m \in [0,\overline{\pi}] $, we set
\begin{align*}
\overline{\mu}_\infty(m) = \E \left[ X_\infty(m) \right].
\end{align*}

\begin{lemma} \label{lemma:mu_bar}
	The function $ m \mapsto \overline{\mu}_\infty(m) $ is continuous, non-decreasing and strictly concave on $ [0,\overline{\pi}] $.
\end{lemma}

\begin{proof}
	Fix $ m_1 \leq m_2 $.
	Then, using the construction in Proposition~\ref{prop:X_inf}, we have, almost surely,
	\begin{align*}
	X_\infty(m_1) \leq X_\infty(m_2).
	\end{align*}
	Taking expectations, we obtain
	\begin{align*}
	\overline{\mu}_\infty(m_1) \leq \overline{\mu}_\infty(m_2).
	\end{align*}
	Hence $ m \mapsto \overline{\mu}_\infty(m) $ is non-decreasing.
	
	The continuity follows from Proposition~\ref{prop:X_inf} and the monotone convergence theorem.
	
	To show that it is concave, we will construct two random variables $ (\nu, X^\delta_\infty(m_1)) $ and $ (\nu, X^\delta_\infty(m_2)) $ distributed according to $ \mu_\infty(m_1 + \delta) $ and $ \mu_\infty(m_2+\delta) $ such that
	\begin{align*}
	X^\delta_\infty(m_1) - X_\infty(m_1) \geq X^\delta_\infty(m_2) - X_\infty(m_2),
	\end{align*}
	almost surely.
	To do this, we will add the same set of global infections (with rate $ \lambda_G \delta \nu / \pib $) to both processes.
	
	Fix $ \delta > 0 $ and let $ \Pi_G^\delta $ be an independent Poisson point process on $ \R_- \times \llbracket 1, \nu \rrbracket $ with intensity $ \delta \frac{\lambda_G}{\pib} dt \otimes c(dk) $.
	We then order the points in $ \Pi_G^\delta $ as above,
	\begin{align*}
	\Pi_G^\delta = \lbrace (t_k^\delta, i_k^\delta) : k \geq 1, 0 > t_1^\delta > t_2^\delta > \ldots \rbrace,
	\end{align*}
	and we define $ I^{k,\delta}(t) $ for $ t \leq 0 $ as above, using the same Poisson point processes of local infections and remission as before, \textit{i.e.} $ \overleftarrow{\Pi}_{rec} $ and $ \overleftarrow{\Pi}_L $.
	We then define
	\begin{align*}
	X^\delta_\infty(m) = \left| \cup_{k \geq 1} \lbrace I^k(0) : u_k \leq m \rbrace \bigcup \cup_{k \geq 1} \lbrace I^{k,\delta}(0) : k \geq 1 \rbrace \right|.
	\end{align*}
	From Proposition~\ref{prop:X_inf}, $ (\nu, X^\delta_\infty(m)) $ is distributed according to $ \mu_\infty(m+\delta) $.
	Furthermore,
	\begin{align*}
	X^\delta_\infty(m) - X_\infty(m) = \left| \cup_{k \geq 1} \lbrace I^{k,\delta}(0), k \geq 1 \rbrace \bigcap \left( \cup_{k \geq 1} \lbrace I^k(0) : u_k \leq m \rbrace \right)^c \right|.
	\end{align*}
	Then, since $ m_1 \leq m_2 $, we have
	\begin{align*}
	\cup_{k \geq 1} \lbrace I^k(0) : u_k \leq m_1 \rbrace \subset \cup_{k \geq 1}
	\lbrace I^k(0) : u_k \leq m_2 \rbrace,
	\end{align*}
	and we deduce that, almost surely,
	\begin{align*}
	X^\delta_\infty(m_1) - X_\infty(m_1) \geq X^\delta_\infty(m_2) - X_\infty(m_2).
	\end{align*}
	Taking expectations, we obtain, for $ m_1 \leq m_2 $,
	\begin{align*}
	\overline{\mu}_\infty(m_1 + \delta) - \overline{\mu}_\infty(m_1) \geq \overline{\mu}_\infty(m_2 + \delta) - \overline{\mu}_\infty(m_2).
	\end{align*}
	This shows that $ m \mapsto \overline{\mu}_\infty(m) $ is concave.
	To show that it is strictly concave, it is sufficient to show that the above inequality is strict with positive probability for any $ \delta > 0 $, which is obvious from our construction.
	This concludes the proof of the lemma.
\end{proof}

\subsection[The basic reproduction number R*]{The basic reproduction number $ R_\ast $} \label{subsec:R0}

Since the non-linear process solves \eqref{eq:X_m} with $ m(t) = \E [X_t(m)] $, if it admits a stationary distribution, we expect that it should be of the form $ \mu_\infty(m) $ with $ m $ satisfying
\begin{align} \label{eq:m_star}
\overline{\mu}_\infty(m) = m.
\end{align}
We note that $ m=0 $ is always a solution to \eqref{eq:m_star}, but, given Lemma~\ref{lemma:mu_bar}, another solution may exist if
\begin{align*}
\frac{d \overline{\mu}_\infty}{dm}(0) > 1.
\end{align*}

\begin{lemma} \label{lemma:R0}
	Recall the definition of $ R_\ast $ in \eqref{def_R0}, then
	\begin{align*}
	\frac{d \overline{\mu}_\infty}{dm}(0) = R_\ast.
	\end{align*}
\end{lemma}

\begin{corollary} \label{cor:m_star}
	If $ R_\ast \leq 1 $, then $ m=0 $ is the unique solution to \eqref{eq:m_star}.
	If $ R_\ast > 1 $, then there exists a unique $ m_\star \in (0,\overline{\pi}] $ satisfying \eqref{eq:m_star}.
\end{corollary}

\begin{proof}
	This is straightforward from Lemma~\ref{lemma:R0} and Lemma~\ref{lemma:mu_bar} and the inequality $ X_\infty(m) \leq \nu $.
\end{proof}

Let us now prove Lemma~\ref{lemma:R0}.

\begin{proof}[Proof of Lemma~\ref{lemma:R0}]
	We prove this result by showing that both terms are equal to the expression given in \eqref{formula_R0}.
	If we set
	\begin{align*}
	\mu^{n,k}_\infty(m) = \P(\nu = n, X_\infty(m) = k),
	\end{align*}
	then the measure $ \mu_\infty(m) $ is characterized by
	\begin{align*}
	\sum_{k=0}^{n} \mathcal{L}\phi(n,k,\pib,m) \mu_\infty^{n,k}(m) = 0, 
	\end{align*}
	for all $\phi:\mathcal{X}\mapsto\mathbb{R}$ measurable and bounded.
	Choosing $ \phi(n,k) = \indic{k \leq \ell} $ for $ 0 \leq \ell \leq n-1 $ yields
	\begin{align*}
	\left( 1 - \frac{\ell}{n} \right) \left[ \lambda_L \frac{n}{n-1}\ell + \lambda_G\frac{n}{\pib} m \right] \mu^{n,\ell}_\infty(m) = \gamma (\ell+1) \mu^{n,\ell+1}_\infty(m).
	\end{align*}
	This, together with the obvious condition $ \sum_{k=0}^{n} \mu^{n,k}_\infty(m) = \pi(n) $ (see \eqref{marginal_pi}) leads to the following expression
	\begin{align*} 
	\begin{split}
	&\mu^{n,\ell}_\infty(m) = \mu^{n,0}_\infty(m) \frac{1}{\gamma^\ell} \prod_{k=0}^{\ell-1} \left\lbrace \frac{1-\frac{k}{n}}{k+1} \left(\lambda_L\frac{n}{n-1} k + \lambda_G \frac{n}{\pib} m\right) \right\rbrace,  \quad 1 \leq \ell \leq n, \\
	&\mu^{n,0}_\infty(m) = \pi(n) \left( 1 + \sum_{\ell=1}^{n} \frac{1}{\gamma^\ell} \prod_{k=0}^{\ell-1} \left\lbrace \frac{1-\frac{k}{n}}{k+1}\left(\lambda_L \frac{n}{n-1}k + \lambda_G \frac{n}{\pib} m\right) \right\rbrace \right)^{-1}.
	\end{split}
	\end{align*}
	From this we deduce easily that
	\begin{align*}
	\frac{d \overline{\mu}_\infty}{dm}(0) = \frac{\lambda_G}{\gamma} \sum_{n=1}^\infty \pi^+(n) \left( 1 + \sum_{\ell=1}^{n-1} \left(\frac{\lambda_L}{\gamma}\right)^{\ell} \prod_{j=1}^\ell \left(1-\frac{j-1}{n-1}\right) \right).
	\end{align*}
	
	We now turn to the quantity $ R_\ast $ defined in \eqref{def_R0}.
	Note that, by the definition of the process $ (I(t), t \geq 0) $ in \eqref{def:I},
	\begin{align} \label{martingale_I}
	\phi(\nu,I(t)) - \phi(\nu,I(0)) - \int_{0}^{t} \mathcal{L}\phi(\nu,I(s),\pib,0)ds
	\end{align}
	is a martingale with respect to the natural filtration of $ \lbrace(\nu,I(t)), t \geq 0\rbrace $.
	Thus if we find a function $ \phi $ such that $ \mathcal{L}\phi(n,x,\pib,0) = x $, we will have
	\begin{align} \label{martingale_I_expect}
	\E\left[ \left. \int_{0}^{T} I(s) ds \right| I(0) = 1, \nu = n \right] = \phi(n,0) - \phi(n,1),
	\end{align}
	where $ T = \inf\lbrace t \geq 0 : I(t) = 0 \rbrace $ (to obtain this, take the expectation of \eqref{martingale_I} at time $ t \wedge T $ and let $ t \to \infty $, using monotone convergence in the integral and dominated convergence in the other term).
	Setting $ \psi(n,x) = \gamma (\phi(n,x-1) - \phi(n,x)) $, $ \mathcal{L}\phi(n,x,\pib,0) = x $ translates into
	\begin{equation*}
	\left\{
	\begin{aligned}
	\psi(n,x)&=1+\frac{\lambda_L}{\gamma}\left(1-\frac{x-1}{n-1}\right)\psi(n,x+1),\ 1\le x\le n-1,\\
	\psi(n,n)&=1.
	\end{aligned}
	\right.
	\end{equation*} 
	We deduce from this that
	\begin{align*}
	\psi(n,1) = \gamma(\phi(n,0)-\phi(n,1)) = 1 + \sum_{\ell=1}^{n-1} \left( \frac{\lambda_L}{\gamma} \right)^{\ell} \prod_{j=1}^\ell \left( 1 - \frac{j-1}{n-1} \right).
	\end{align*}
	Together with \eqref{martingale_I_expect}, this proves the lemma.
\end{proof}

Let us quickly mention another avenue for proving Lemma~\ref{lemma:R0}, which makes use of Proposition~\ref{prop:X_inf}.
For $ \varepsilon > 0 $, let us write
\begin{align*}
\lbrace k \geq 1 : u_k \leq \varepsilon \rbrace = \lbrace 1 \leq k_1(\varepsilon) < k_2(\varepsilon) < \ldots \rbrace.
\end{align*}
Then we write
\begin{align*}
X_\infty(\varepsilon) = | I^{k_1(\varepsilon)}(0) | + | \cup_{j \geq 2} I^{k_j(\varepsilon)}(0) \cap I^{k_1(\varepsilon)}(0)^c |.
\end{align*}
Then, noting that $ -t_{k_1(\varepsilon)} $ is distributed as an exponential variable with parameter $ \lambda_G\nu  \varepsilon/\pib $, it is possible to show that
\begin{align*}
\E\left[ \left. | I^{k_1(\varepsilon)}(0) | \right| \nu = n \right] = \varepsilon \lambda_G\frac{n}{\pib} \int_{0}^{\infty} \E\left[ \left. | I^{1}(t_1 + t) | \right| \nu = n \right] dt + \mathit{o}(\varepsilon),
\end{align*}
and that
\begin{align*}
\E\left[ \left. | \cup_{j \geq 2} I^{k_j(\varepsilon)}(0) \cap I^{k_1(\varepsilon)}(0)^c | \right| \nu = n \right] = \mathit{o}(\varepsilon).
\end{align*}
We then finish by noting that $ I^1(t_1+t) $ is distributed as $ I(t) $ conditionally on $ I(0) = 1 $ and that
\begin{align*}
\frac{d \overline{\mu}_\infty}{dm}(0) = \lim_{\varepsilon \downarrow 0} \frac{1}{\varepsilon} \E[X_\infty(\varepsilon)].
\end{align*}

\subsection{Large-time behaviour of the non-linear Markov process}

We now prove Theorem~\ref{thm:large_time}.
We split the proof in two parts, first dealing with the case $ R_\ast \leq 1 $ and then with $ R_\ast > 1 $.

\begin{proof}[Proof of Theorem~\ref{thm:large_time}, $ R_\ast \leq 1 $]
	Let $ m^+_0(t) = \overline{\pi} $ and set, for $ k \geq 0 $,
	\begin{align*}
	m^+_{k+1}(t) = \overline{\mu}_t(m^+_k,\mu_0).
	\end{align*}
	Clearly $ \E[X(t)] \leq m^+_0(t) $ for all $ t \geq 0 $.
	Since $ (\E[X(t)], t \geq 0) $ is a fixed point of $ m(\cdot) \mapsto \overline{\mu}_\cdot(m,\mu_0) $ and using Lemma~\ref{lemma:monotonicity}, for every $ k \geq 0 $,
	\begin{align} \label{encadrement_1}
	0 \leq \E[X(t)] \leq m^+_k(t).
	\end{align}
	Furthermore, by Lemma~\ref{lemma:limit_forced_process}, for all $ k \geq 0 $,
	\begin{align*}
	\lim_{t \to \infty} m^+_k(t) = \overline{\mu}_\infty^{\circ k}(\overline{\pi}),
	\end{align*}
	where $ \overline{\mu}^{\circ k}_\infty(\cdot) = \overline{\mu}_\infty(\overline{\mu}_\infty(\ldots)) $ is the $ k $-th iterate of $ m \mapsto \overline{\mu}_\infty(m) $.
	Letting $ t \to \infty $ in \eqref{encadrement_1},
	\begin{align*}
	0 \leq \liminf_{t \to \infty} \E[X(t)] \leq \limsup_{t \to \infty} \E[X(t)] \leq \overline{\mu}_\infty^{\circ k}(\overline{\pi}).
	\end{align*}
	But, by Lemma~\ref{lemma:mu_bar} and Lemma~\ref{lemma:R0}, since $ R_\ast \leq 1 $,
	\begin{align*}
	\overline{\mu}_\infty^{\circ k}(\overline{\pi}) \to 0 \text{ as } k \to \infty.
	\end{align*}
	As a result,
	\begin{align*}
	\lim_{t \to \infty} \E[X(t)] = 0,
	\end{align*}
	and the result follows.
\end{proof}

Before proving the result when $ R_\ast > 1 $, we state the following lemma, whose proof we delay until Subsection~\ref{subsec:liminf}.

\begin{lemma} \label{lemma:liminf}
	Suppose that $ R_\ast > 1 $ and that $ \E[X(0)] > 0 $, then
	\begin{align*}
	\liminf_{t \to \infty} \E[X(t)] > 0.
	\end{align*}
\end{lemma}

Let us now finish the proof of Theorem~\ref{thm:large_time}.

\begin{proof}[Proof of Theorem~\ref{thm:large_time}, $ R_\ast > 1 $]
	The strategy of the proof is similar to the case $ R_\ast \leq 1 $, but we now define two functions
	\begin{align*}
	m^+_0(t) = \overline{\pi}, && m^-_0(t) = \inf_{s \geq 0}\E[X(s)].
	\end{align*}
	Note that by Lemma~\ref{lemma:liminf}, $ \lim_{t \to \infty} m_0^-(t) > 0 $.
	As before, we set, for $ k \geq 0 $,
	\begin{align*}
	m^+_{k+1}(t) = \overline{\mu}_t(m^+_k,\mu_0), && m^-_{k+1}(t) =\overline{\mu}_t(m^-_k,\mu_0).
	\end{align*}
	By the same argument as before, since $ m^-_0(t) \leq \E[X(t)] \leq m^+_0(t) $, we have, for every $ k \geq 0 $,
	\begin{align*}
	m^-_{k}(t) \leq \E[X(t)] \leq m^+_k(t).
	\end{align*}
	Using Lemma~\ref{lemma:limit_forced_process} and letting $ t \to \infty $, we obtain
	\begin{align} \label{encadrement_2}
	\overline{\mu}_\infty^{\circ k}(\inf_{t \geq 0}\E[X(t)]) \leq \liminf_{t \to \infty} \E[X(t)] \leq \limsup_{t \to \infty}\E[X(t)] \leq \overline{\mu}_\infty^{\circ k}(\overline{\pi}).
	\end{align}
	But, by Lemma~\ref{lemma:mu_bar} and the fact that $ R_\ast > 1 $, we have
	\begin{align*}
	\lim_{k \to \infty} \overline{\mu}_\infty^{\circ k}(\inf_{t \geq 0}\E[X(t)]) = \lim_{k \to \infty} \overline{\mu}_\infty^{\circ k}(\overline{\pi}) = m_\star,
	\end{align*}
	where $ m_\star \in (0,\overline{\pi}] $ is defined by Corollary~\ref{cor:m_star} (also using Lemma~\ref{lemma:liminf} and the fact that $ \E[X(0)] > 0 $).
	Hence, letting $ k \to \infty $ in \eqref{encadrement_2},
	\begin{align*}
	\E[X(t)] \to m_\star,
	\end{align*}
	as $ t \to \infty $.
	Finally by Lemma~\ref{lemma:limit_forced_process}, since the non-linear process is the forced process with $ m(t) = \E[X(t)] $,
	\begin{align*}
	(\nu, X(t)) \to \mu_\infty(m_\star),
	\end{align*}
	in distribution as $ t \to \infty $, and the theorem is proved.
\end{proof}

Note that, without Lemma~\ref{lemma:liminf}, we would not have been able to bound $ \E[X(t)] $ from below by anything useful, since $ \overline{\mu}_\infty(0) = 0 $.

\subsection{Branching process minoration} \label{subsec:liminf}

\begin{proof}[Proof of Lemma~\ref{lemma:liminf}]
	For $ c > 0 $ and $ N_0 \geq 1 $, let $ \pi_c(n) $ be defined by
	\begin{align*}
		\pi_c(n) :=  \frac{{\bf 1}_{n \leq N_0} e^{-cn} \pi^+(n)}{\sum_{m=1}^{N_0} e^{-cm} \pi^+(m)}, \quad n \geq 1.
	\end{align*}
	Then $ \pi_c $ stochastically dominates $ \pi_{c'} $ as soon as $ c \leq c' $ and $ \pi^+ $ stochastically dominates $ \pi_c $ for any $ c > 0 $.
	Since $ R_\ast > 1 $, we can choose $ N_0 \geq 1 $, $ p \in (0,1) $ and $ c > 0 $ such that
	\begin{align} \label{reprod_nb_bp}
	(1-p) \lambda_G \sum_{n=1}^{N_0} \pi_c(n) \E\left[ \left. \int_{0}^{\infty} I(t) dt \right| \nu=n, I(0) = 1 \right] > 1.
	\end{align}
	We then choose $ q \in (0,p) $ such that
	\begin{align} \label{q_1}
		q < \frac{\pib}{\pibp} p^2,
	\end{align}
	and, unless $ \pi $ is concentrated on a single point,
	\begin{align} \label{q_2}
		q < \frac{\pib}{\pibp} \left( \frac{d_c}{1+d_c} \right)^2,
	\end{align}
	where
	\begin{align*}
		d_c = \inf \left\lbrace \sum_{\ell = 1}^{n} \pi_c(\ell) - \sum_{\ell=1}^{n} \pi^+(\ell), 1 \leq n < N_0 : \pi(n) > 0 \right\rbrace.
	\end{align*}
	Since $ \pi^+ $ stochastically dominates $ \pi_c $, $ d_c > 0 $ as soon as $ \pi $ is not concentrated on a single point (if $ \pi $ has compact support, then $ N_0 $ should not be larger than the upper bound of this support). 
	
	Without loss of generality, we can assume that there exist $ N_1, N_2 $ in $ \N $ such that $ p = 1/N_1 $ and $ q = 1/N_2 $.
	For the rest of this proof, we restrict $ N $ to multiples of both $ N_1 $ and $ N_2 $.
	
	Let $Z^N_t$ denote the number of households containing at least one infectious individual at time $ t $, \textit{i.e.}
	\begin{align*}
	Z^N_t=\sum_{i=1}^N{\bf1}_{X^N_i(t)\ge1},
	\end{align*}
	where $\{X^N_i(t),\ t\ge0;\ 1\le i\le N\}$ is the solution of equation \eqref{Nsde}.
	
	We now define a continuous-time non-Markovian branching process of household infections as follows.
	Start with $ Y^N_0 = Nq $ infected households, each with a single infectious individual.
	If there are currently $ k $ infected households with $ x_1, \ldots, x_k $ infectious individuals, at rate $ (1-p) \lambda_G \sum_{i=1}^k x_i $, a new household, whose size is chosen according to the distribution $ \pi_c $, is added to the process with a single infectious individual.
	Apart from this, each household undergoes a local epidemic with rates $ \lambda_L $ and $ \gamma $, independently from the others.
	Then $ Y^N_t $ denotes the number of infected households at time $ t \geq 0 $.
	
	The corresponding discrete time branching process is supercritical, since the expected number of ``offspring'' of each household is given by \eqref{reprod_nb_bp}, which has been chosen to be larger than one. 
	Then from Lemma 2.1 in Doney \cite{Do}, if $r>0$ denotes the real number such that 
	\begin{align*}
	(1-p) \lambda_G \sum_{n=1}^{N_0} \pi_c(n) \int_0^\infty e^{-r t} \E_1[I(t) | \nu = n]dt=1,
	\end{align*}
	where $(I(t), t \geq 0)$ is the process defined in \eqref{def:I} and 
	$\E_1$ means that we take the expectation under the initial condition $I(0)=1$, then there exists $ a > 0 $ such that
	\begin{align} \label{asymptotic_Y}
	\E[Y^N_t] \sim N a e^{rt} \text{ as } t \to \infty.
	\end{align}
	
	Suppose that $ Nq\le Z^N_0 $, and that the initially infected households in $ Y^N_0 $ have the same sizes as those in $ Z^N_0 $. 
	Then define
	\begin{align*}
	T_{N,p} = \inf \left\lbrace t \geq 0 : \frac{\sum_{i=1}^{N} \nu_i {\bf1}_{X^N_i(t) \geq 1}}{\sum_{i=1}^{N} \nu_i} > p \right\rbrace
	\end{align*}
	and
	\begin{align*}
		\tau_{N,c} = \inf \left\lbrace t \geq 0 : \exists n \geq 1, \frac{\sum_{i=1}^{N} \nu_i {\bf1}_{\nu_i \leq n, X^N_i(t) = 0}}{\sum_{i=1}^{N} \nu_i {\bf1}_{X^N_i(t) = 0}}  > \sum_{\ell = 1}^{n} \pi_c(\ell) \text{ or } \sum_{i=1}^{N} {\bf1}_{X^N_i(t) = 0} = 0 \right\rbrace.
	\end{align*}
	We claim that, on the interval $[0,T_{N,p} \wedge \tau_{N,c})$, $Z^N_t$ stochastically dominates $Y^N_t$ (\textit{i.e.} we can define $ (Y^N_t, t \geq 0) $ such that $ Y^N_t \leq Z^N_t $ for $ t \in [0,T_{N,p}\wedge\tau_{N,c}) $).
	To see this, note that $ Y^N_0 \leq Z^N_0 $ and that, since each household in $ Y^N_t $ starts with a single infectious individual, the number of infectious individuals in each household can be kept larger in $ Z^N_0 $ than in $ Y^N_0 $.
	This stays true until the first time at which a \textit{new} household is infected in either process, since the local infection parameters are the same in both processes, and in $ Z^N $, there are additional infections due to global infections between already infected households.
	Furthermore, in the process $ (Z^N_t, t \geq 0) $, a \textit{new} household is infected at rate
	\begin{align*}
	\lambda_G \frac{1}{N} \sum_{j=1}^{N} X^N_j(t) \sum_{i=1}^{N} \frac{\nu_i}{\overline{\nu}^N} {\bf1}_{X^N_i(t) = 0} = \lambda_G \left( 1 - \frac{\sum_{i=1}^{N} \nu_i {\bf1}_{X_i^N(t) \geq 1}}{\sum_{i=1}^{N} \nu_i} \right) \sum_{j=1}^{N} X^N_j(t).
	\end{align*}
	and for $ t \in [0,T_{N,p}) $, this rate is larger than the rate at which a new household is infected in the process $ (Y^N_t, t \geq 0) $, provided the total number of infectious individuals is larger in the former process.
	To achieve this, it is enough to be able to choose the sizes of the new households so that, each time a new household is added to the branching process, the household that is added to the original process at the same time is larger than the first (two households undergoing local infections at the same rate can be coupled so that the largest one always contains more infected individuals than the other).
	Since the newly infected households in the original process are chosen proportionally to their size, while the new households in the branching process are chosen according to the distribution $ \pi_c $, this is possible for all $ t \leq \tau_{N,c} $.
	We can thus couple the two processes in such a way that 
	\begin{align*}
	Y^N_t \leq Z^N_t, \quad \forall t \in [0,T_{N,p} \wedge \tau_{N,c}),
	\end{align*}
	almost surely for all $N\ge1$. 
	
	Now, by Theorem~\ref{thm:chaos}, as $ N \to \infty $, for any $ T > 0 $,
	\begin{align} \label{cvg_ZN}
	\frac{Z^N_t}{N} \to \mathfrak{p}(t):=\P(X(t)\ge1),
	\end{align}
	and
	\begin{align} \label{cvg_rate}
	\frac{\sum_{i=1}^{N} \nu_i {\bf1}_{X_i^N(t) \geq 1}}{\sum_{i=1}^{N} \nu_i} \to \frac{1}{\pib} \E[\nu {\bf1}_{X(t) \geq 1}], && \frac{\sum_{i=1}^{N} \nu_i {\bf1}_{\nu_i \leq n, X^N_i(t) = 0}}{\sum_{i=1}^{N} \nu_i {\bf1}_{X^N_i(t) = 0}} \to \frac{\E\left[ \nu {\bf 1}_{\nu \leq n, X(t) = 0} \right]}{\E\left[ \nu {\bf 1}_{X(t) = 0} \right]},
	\end{align}
	uniformly on $ [0,T] $, in probability.
	Furthermore, there exists a deterministic function $ f : \R_+ \to \R_+ $ such that
	\begin{align*}
	\frac{Y^N_t}{N} \to f(t),
	\end{align*}
	uniformly on $ [0,T] $ as $ N \to \infty $, in probability.
	Furthermore, by \eqref{asymptotic_Y}, 
	\begin{align} \label{asymptotic_f}
	f(t) \sim a e^{rt} \text{ as } t \to \infty.
	\end{align}
	
	For any $ p' < p $, let
	\begin{align*}
	T_{p'} = \inf \left\lbrace t \geq 0 : \frac{1}{\pib}\E[\nu {\bf1}_{X(t) \geq 1}] > p' \right\rbrace.
	\end{align*}
	Also define, for $ c' \in (0,c) $,
	\begin{align*}
		\tau_{c'} = \inf \left\lbrace t \geq 0 : \exists n \geq 1, \frac{\E\left[ \nu {\bf 1}_{\nu \leq n, X(t) = 0} \right]}{\E\left[ \nu {\bf 1}_{X(t) = 0} \right]} > \sum_{\ell=1}^{n} \pi_{c'}(\ell) \right\rbrace.
	\end{align*}
	By \eqref{cvg_rate} and the fact that, for $ c' < c $, $ \pi_{c'} $ stochastically dominates $ \pi_c $, for any $ t \leq T_{p'} \wedge \tau_{c'} $, $t<\liminf_{N\to\infty}T_{N,p} \wedge \tau_{N,c}$, and consequently for $N$ large enough,
	\begin{align*}
	\frac{Z^N_t}{N} \geq \frac{Y^N_t}{N},\ \forall t\le T_{p'} \wedge \tau_{c'} \wedge T,
	\end{align*}
	for any $ T > 0 $ (noting that the probability that all households are simultaneously infected for some time $ t $ vanishes as $ N \to \infty $ since we have assumed $ Z^N_0 = pN $).
	Letting $ N \to \infty $, we obtain
	\begin{align*}
	\mathfrak{p}(t) \geq f(t),\ \forall t \leq T_{p'} \wedge \tau_{c'}.
	\end{align*}
	Now define, for some $ b > q $,
	\begin{align*}
	T^f_b = \inf \lbrace t \geq 0 : f(t) > b \rbrace.
	\end{align*}
	Then, if $ T^f_b < T_{p'} \wedge \tau_{c'} $, $ \mathfrak p(T^f_b) \geq b $.
	If however $ T_{p'} \leq T^f_b $, then, by the Cauchy-Schwarz inequality,
	\begin{align} \label{cauchy-schwarz}
	\E[\nu {\bf1}_{X(t) \geq 1}] \leq \sqrt{\E[\nu^2 ]} \sqrt{\mathfrak{p}(t)},
	\end{align}
	and thus,
	\begin{align*}
	\mathfrak p(T_{p'}) \geq (p')^2 \frac{\pib}{\pibp}.
	\end{align*}
	By \eqref{q_1}, we can choose $ p' $ such that the right hand side is larger than $ q $.
	Finally, we note that
	\begin{align*}
	\frac{\mathbb{E}\left[ \nu {\bf 1}_{\nu = \ell, X(t) = 0} \right]}{\mathbb{E}\left[ \nu {\bf 1}_{X(t) = 0} \right]} - \pi^+(\ell) &= \frac{\mathbb{E}\left[ \nu {\bf 1}_{\nu = \ell, X(t) = 0} \right] - \mathbb{E}\left[ \nu {\bf 1}_{\nu = \ell} \right]}{\mathbb{E}\left[ \nu {\bf 1}_{X(t) = 0} \right]} + \ell \pi(\ell) \frac{\mathbb{E}\left[\nu\right] - \mathbb{E}\left[ \nu {\bf 1}_{X(t) = 0} \right]}{\overline{\pi} \mathbb{E}\left[ \nu {\bf 1}_{X(t) = 0} \right]} \\
	&= - \frac{\mathbb{E}\left[ \nu {\bf 1}_{\nu = \ell, X(t) \geq 1} \right]}{\mathbb{E}\left[ \nu {\bf 1}_{X(t) = 0} \right]} + \pi^+(\ell) \frac{\mathbb{E}\left[ \nu {\bf 1}_{X(t) \geq 1} \right]}{\mathbb{E}\left[ \nu {\bf 1}_{X(t) = 0} \right]} \\
	&\leq \pi^+(\ell) \frac{\mathbb{E}\left[ \nu {\bf 1}_{X(t) \geq 1} \right]}{\mathbb{E}\left[\nu {\bf 1}_{X(t) = 0} \right]}.
	\end{align*}
	Hence if $ \tau_{c'} < T^f_b $, there exists $ 1 \leq n < N_0 $ such that $ \pi(n) > 0 $ and
	\begin{align*}
	d_{c'} \leq \sum_{\ell = 1}^{n} \pi_{c'}(\ell) - \sum_{\ell = 1}^{n} \pi^+(\ell) \leq \frac{\mathbb{E}\left[ \nu {\bf 1}_{X(\tau_{c'}) \geq 1} \right]}{\mathbb{E}\left[\nu {\bf 1}_{X(\tau_{c'}) = 0}\right]} \sum_{\ell = 1}^{n} \pi^+(\ell) \leq \frac{\mathbb{E}\left[ \nu {\bf 1}_{X(\tau_{c'}) \geq 1} \right]}{\mathbb{E}\left[\nu {\bf 1}_{X(\tau_{c'}) = 0}\right]}.
	\end{align*}
	Since $ \mathbb{E}\left[\nu {\bf 1}_{X(\tau_{c'}) = 0}\right] = \overline{\pi} - \mathbb{E}\left[ \nu {\bf 1}_{X(\tau_{c'}) \geq 1} \right] $, this yields
	\begin{align*}
		\mathbb{E}\left[ \nu {\bf 1}_{X(\tau_{c'}) \geq 1} \right] \geq \overline{\pi} \frac{d_{c'}}{1+d_{c'}}.
	\end{align*}
	Thus, using \eqref{cauchy-schwarz} again,
	\begin{align*}
	\mathfrak p(\tau_{c'}) \geq \frac{\overline{\pi}}{\overline{\pi}^+} \left(\frac{d_{c'}}{1+d_{c'}}\right)^2,
	\end{align*}
	where $ c' $ can be chosen so that the right hand side is larger than $ q $, by \eqref{q_2} (note that, if $ \pi $ is concentrated on a single point, then $ \tau_{c'} = +\infty $).
	
	As a consequence, if for some $t\ge0$, $\mathfrak{p}(t)=q$, then $\mathfrak p(t+s)$ reaches $ b \wedge (p')^2 \overline{\pi} / \overline{\pi}^+ \wedge (d_{c'}/(1+d_{c'}))^2 \overline{\pi}/\overline{\pi}^+ > q $ for some $ s \leq T^f_b $, and $ \mathfrak p(t+u) \geq f(u) $ for all $ u \in [0,s] $.
	By \eqref{asymptotic_f}, $ f $ is uniformly bounded away from 0 on $ [0,T^f_b] $.
	This proves the Lemma.
\end{proof}

\begin{acks}[Acknowledgments]
The authors wish to thank Frank Ball for pointing out to them that their first version of this work was not consistent with the classical household models studied in the literature, and two anonymous referees, whose comments resulted in several improvements in our paper.
\end{acks}

\begin{funding}
The two authors were supported in part by the chair “Mod\'elisation math\'ematique et Biodiversit\'e” of Veolia-Ecole Polytechnique-Mus\'eum National d’Histoire Naturelle-Fondation X.
\end{funding}

\end{document}